# Discrete q-Hermite polynomials: An elementary approach

Johann Cigler

Fakultät für Mathematik, Universität Wien
Uni Wien Rossau, Oskar-Morgenstern-Platz 1, 1090 Wien
johann.cigler@univie.ac.at

**Abstract**

We present a simple approach to discrete $q-$ Hermite polynomials with special emphasis on analogies with the classical case.

## 0. Introduction

The so-called "discrete $q-$ Hermite polynomials" are $q-$ analogues of the classical Hermite polynomials which generalize most of their elementary properties. The purpose of this paper is to emphasize these analogies without recourse to sophisticated general theories. Some 35 years ago I stumbled upon these polynomials (cf. [3], [4]) as I studied some $q-$ identities from the point of view of Rota's umbral calculus ([15]). In his thesis Peter Kirschenhofer [12] obtained further properties and determined a measure on the real line with respect to which they are orthogonal. At that time I was unaware that Al-Salam and Carlitz [1] had introduced these polynomials some years earlier via generating functions. Some properties of these polynomials are collected in [13] within the framework of the Askey- scheme of hypergeometric orthogonal polynomials and its $q-$ analogue. They are also dealt with in an elementary way in [8], [9] and [14]. Here I give an introduction to these polynomials in the spirit of my papers [3], [4]. I shall derive recurrence relations, generating functions, Rodrigues type formulae together with $q-$ analogues of Burchnall's and Nielsen's formula and some measures with respect to which the polynomials are orthogonal. Finally a curious connection with tangent and Euler numbers is mentioned. Whereas most results are well known or hidden in more general theories the approach seems to be novel.

## 1. Background material about Hermite polynomials

In this section I state some well-known results about the classical (probabilists') Hermite polynomials (cf. [2], [15]) which have beautiful generalizations for discrete $q-$ Hermite polynomials. The main interest is in $H_n(x)=H_n(x,1)$ but it is convenient to state some results more generally for bivariate Hermite polynomials.

For $s>0$ the bivariate Hermite polynomials

$$H_n(x,s)=\sum_{k=0}^{\left\lfloor \frac{n}{2}\right\rfloor}(-1)^k\binom{n}{2k}(2k-1)!!s^k x^{n-2k} \tag{1.1}$$

are monic polynomials which satisfy

$$\frac{d}{dx}H_n(x,s)=nH_{n-1}(x,s) \tag{1.2}$$

and are orthogonal with respect to the linear functional $\Lambda$ defined by

$$\Lambda\left(H_n(x,s)\right)=[n=0]. \tag{1.3}$$

They satisfy the recurrence

$$H_{n+1}(x,s)=xH_n(x,s)-nsH_{n-1}(x,s) \tag{1.4}$$

with initial values $H_0(x,s)=1$ and

Their generating function is

$$\sum_{n\ge 0}H_n(x,s)\frac{z^n}{n!}=e^{xz}e^{-\frac{sz^2}{2}}. \tag{1.5}$$

Since the differentiation operator $D$ satisfies $De^{xz}=ze^{xz}$ this can also be written as $\sum_{n\ge 0}H_n(x,s)\frac{z^n}{n!}=e^{-\frac{sD^2}{2}}e^{xz}$ or equivalently

$$H_n(x,s)=e^{-\frac{sD^2}{2}}x^n. \tag{1.6}$$

By (1.5) we have $\Lambda\left(e^{xz}\right)=e^{\frac{sz^2}{2}}$ which by comparing coefficients is equivalent with the moments

$$\begin{aligned}\Lambda\left(x^{2m}\right)&=s^m(2m-1)!!\\ \Lambda\left(x^{2m+1}\right)&=0.\end{aligned} \tag{1.7}$$

The Taylor expansion

$$\sum_{n\ge 0}H_n(x,s)\frac{z^n}{n!}=e^{xz}e^{-\frac{sz^2}{2}}=e^{\frac{x^2}{2s}}e^{-\frac{(x-sz)^2}{2s}}=e^{\frac{x^2}{2s}}\sum_{n\ge 0}\left(D^n e^{\frac{-x^2}{2s}}\right)\frac{(-sz)^n}{n!}$$

gives the Rodrigues formula

$$H_n(x,s)=\frac{1}{e^{-\frac{x^2}{2s}}}\left(-sD\right)^n e^{-\frac{x^2}{2s}}. \tag{1.8}$$

Let **g(x)** denote the multiplication operator with a function $g(x)$, i.e. **g(x)** $f(x)=g(x)f(x)$.

Since $D\mathbf{g}(\mathbf{x}) = \mathbf{g}(\mathbf{x})D + \mathbf{g}'(\mathbf{x})$ we have $sD\mathbf{e}^{-\frac{\mathbf{x}^2}{2\mathbf{s}}} = \mathbf{e}^{-\frac{\mathbf{x}^2}{2\mathbf{s}}}sD - x\mathbf{e}^{-\frac{\mathbf{x}^2}{2\mathbf{s}}}$ and therefore

$$\frac{\mathbf{1}}{\mathbf{e}^{-\frac{\mathbf{x}^2}{2\mathbf{s}}}}\left(-sD\right)\mathbf{e}^{-\frac{\mathbf{x}^2}{2\mathbf{s}}} = \mathbf{x} - sD. \tag{1.9}$$

This implies that the Hermite polynomials $H_n(x,s)$ satisfy

$$H_n(x,s) = \left(\mathbf{x} - sD\right)^n 1. \tag{1.10}$$

A generalization of (1.10) is Burchnall's operational formula (cf. e.g. [15])

$$\left(\mathbf{x} - sD\right)^n = \sum_{k=0}^{n} \binom{n}{k} \mathbf{H}_{\mathbf{n}\text{-}\mathbf{k}}\left(\mathbf{x}, s\right)\left(-sD\right)^k, \tag{1.11}$$

which implies Nielsen's identity

$$H_{n+m}(x,s) = \sum_k \binom{n}{k}\binom{m}{k} k!(-s)^k H_{n-k}(x,s)H_{m-k}(x,s). \tag{1.12}$$

The Rodrigues formula provides an easy proof that

$$\Lambda\left(f\right) = \frac{1}{\sqrt{2\pi s}}\int_{-\infty}^{\infty} f(x)e^{-\frac{x^2}{2s}}dx \tag{1.13}$$

for polynomials $f(x)$.

It suffices to show that $\int_{-\infty}^{\infty} H_n(x,s)e^{-\frac{x^2}{2s}}dx = 0$ for $n > 0$.

This is immediate by (1.8) and integration by parts:

$$\int_{-\infty}^{\infty} H_n(x,s)e^{-\frac{x^2}{2s}}dx = (-s)^n\int_{-\infty}^{\infty} D^n e^{-\frac{x^2}{2s}}dx = -(-s)^n\, D^{n-1}e^{-\frac{x^2}{2s}}\Big|_{-\infty}^{\infty} = 0.$$

For $n = 0$ we need the normalization $\frac{1}{\sqrt{2\pi s}}\int_{-\infty}^{\infty} e^{-\frac{x^2}{2s}}dx = 1.$

From (1.4) we see that

$$\Lambda\left(H_n(x,s)H_m(x,s)\right) = s^n n![n = m]. \tag{1.14}$$

## 2. Bivariate q-Hermite polynomials and their generating functions

### 2.1. The polynomials

In this section I freely use some well-known notations of $q-$analysis in the form I used in [7]. Let me only mention that we always assume $0<q<1$ and that $D_q$ denotes the $q-$differentiation operator defined by $D_q f(\mathrm{x})=\frac{f(x)-f(qx)}{(1-q)x}$ which satisfies $D_q x^n=[n]x^{n-1}$ with $[n]=[n]_q=\frac{1-q^n}{1-q}.$

My original approach to these polynomials has been (cf. [3], [4]) to look for sequences of monic polynomials which satisfy

$$D_q p_n(x)=[n]p_{n-1}(x) \tag{2.1}$$

and are orthogonal with respect to a linear functional $\Lambda.$ This means that $\Lambda\left(p_m(x)p_n(x)\right)=c_n[n=m]$ with $c_n\neq 0.$

By Favard's theorem they must satisfy a 3-term recurrence of the form

$$p_{n+1}(x)=xp_n(x)+a(n)p_{n-1}(x). \tag{2.2}$$

If we apply the operator $D_q$ and observe that $D_q\left(xf(x)\right)=qxD_q f(x)+f(x)$ we get

$$[n+1]p_n(x)=qx[n]p_{n-1}(x)+p_n(x)+a(n)\left[n-1\right]p_{n-2}(x)$$

or

$$p_n(x)=xp_{n-1}(x)+\frac{a(n)\left[n-1\right]}{q[n]}p_{n-2}(x).$$

Comparing with (2.2) we see that

$$\frac{a(n)\left[n-1\right]}{q[n]}=a(n-1)$$

or

$$a(n)=q^{n-1}[n]a(1).$$

By choosing $a(1)=-qs$ we get the $q-$ Hermite polynomials $H_n(x,s,q)$ which satisfy

$$H_{n+1}(x,s,q)=xH_n(x,s,q)-q^n s[n]H_{n-1}(x,s,q) \tag{2.3}$$

with initial values $H_{-1}(x,s,q)=0$ and $H_0(x,s,q)=1.$

It is easy to determine the coefficients of $H_n(x,s,q).$

From (2.3) it is clear that $H_n(x,s,q)=\sum a(n,k)x^{n-2k}$ for some $a(n,k).$

Futhermore $H_n(0,s,q)=-q^{n-1}s[n-1]H_{n-2}(0,s,q)$ which gives $H_{2k+1}(0,s,q)=0$ and $H_{2k}\left(0,s,q\right)=(-1)^k s^k q^{k^2}[2k-1]!!.$

This implies $[n-2k]!a[n,k]= D_q^{n-2k}\left(H_n(x,s,q)\right)\Big|_{x=0}=\frac{[n]!}{[2k]!}H_{2k}(0,s,q)$ and therefore

$$H_n(x,s,q)=\sum_{2k\le n}(-s)^k q^{k^2}\begin{bmatrix} n \\ 2k \end{bmatrix}[2k-1]!!x^{n-2k}. \tag{2.4}$$

It is easily verified that (2.4) really satisfies both (2.1) and (2.3).

Now we define the linear functional $\Lambda$ by $\Lambda\left(H_n(x,s,q)\right)=[n=0].$ By Favard's theorem we know that the sequence $\left(H_n(x,s,q)\right)$ is orthogonal with respect to $\Lambda.$

Since $H_n(x,s,q)$ is monic we have $\Lambda\left(H_n(x,s,q)^2\right)=\Lambda\left(x^n H_n(x,s,q)\right).$

By (2.3) we see that $\Lambda\left(x^n H_n(x,s,q)\right)=q^n s[n]\Lambda\left(x^{n-1}H_{n-1}(x,s,q)\right).$

This implies

$$\Lambda\left(H_n(x,s,q)H_m(x,s,q)\right)=q^{\binom{n+1}{2}}s^n[n]![n=m]. \tag{2.5}$$

From (2.3) we deduce the following determinant representation

$$H_n(x,s,q)=\det\begin{pmatrix} x & qs & 0 & \cdots & 0 & 0 \\ 1 & x & [2]q^2 s & \cdots & 0 & 0 \\ 0 & 1 & x & \cdots & 0 & 0 \\ \vdots & \vdots & \vdots & \ddots & \vdots & \\ 0 & 0 & 0 & \cdots & x & [n-1]q^{n-1}s \\ 0 & 0 & 0 & \cdots & 1 & x \end{pmatrix}$$

If we replace $q$ by $\frac{1}{q}$ and observe that $\left[n\right]_{\frac{1}{q}}=q^{1-n}\left[n\right]_q$ we get

$$H_{n+1}(x,s,q^{-1})=xH_n(x,s,q^{-1})-q^{1-2n}s[n]_q H_{n-1}(x,s,q^{-1}) \tag{2.6}$$

and

$$H_n(x,s,q^{-1})=\sum_{2k\le n}(-s)^k q^{\binom{n-2k}{2}-\binom{n}{2}}\begin{bmatrix} n \\ 2k \end{bmatrix}[2k-1]!!x^{n-2k}. \tag{2.7}$$

Since $D_{q^{-1}}=\varepsilon^{-1}D_q$ identity (2.1) becomes

$$D_q H_n\left(x,s,\frac{1}{q}\right)=\frac{1}{q^{n-1}}[n]_q H_{n-1}\left(qx,s,\frac{1}{q}\right). \tag{2.8}$$

Some formulae become simpler if we consider the polynomials

$$K_n(x,s,q)=q^{\binom{n}{2}}H_n\left(x,s,\frac{1}{q}\right). \tag{2.9}$$

Here we have

$$K_n\left(x,s,q\right)=\sum_k(-s)^k\begin{bmatrix}n\\2k\end{bmatrix}[2k-1]!!q^{\binom{n-2k}{2}}x^{n-2k} \tag{2.10}$$

and

$$K_{n+1}(x,s,q)=q^n xK_n(x,s,q)-s[n]K_{n-1}(x,s,q). \tag{2.11}$$

Furthermore

$$\left(D_q K_n\right)(x,s,q)=[n]K_{n-1}(qx,s,q). \tag{2.12}$$

**Remark**

There are many orthogonal $q-$ polynomials which reduce to $H_n(x,s)$ for $q\to 1.$ But only the continuous $q-$ Hermite polynomials $\tilde{H}_n(x,s,q)$ (cf. e.g. [7] ) which satisfy $\tilde{H}_{n+1}(x,s,q)=x\tilde{H}_n(x,s,q)-[n]s\tilde{H}_{n-1}(x,s,q)$ seem to have interesting properties too.

On the other hand there are $q-$ analogues with simple formulae and recurrence relations but which are not orthogonal. For example the simplest $q-$ analogue of formula (1.1) seems to be

$$\bar{h}_n(x,s,q)=\sum_{k=0}^{\left\lfloor\frac{n}{2}\right\rfloor}\begin{bmatrix}n\\2k\end{bmatrix}\left[2k-1\right]!!(-s)^k x^{n-2k} \tag{2.13}$$

which has been studied by Kirschenhofer [12] and satisfies

$$D_q\bar{h}_n(x,s,q)=[n]\bar{h}_n(x,s,q), \tag{2.14}$$

$$\bar{h}_{n+1}(x,s,q)=x\bar{h}_n(x,s,q)-s[n]\bar{h}_{n-1}(x,s,q)+xs[n](1-q^{n-1})\bar{h}_{n-2}(x,s,q) \tag{2.15}$$

and

$$\bar{h}_{n+1}(x,s,q)=x\bar{h}_n(x,s,q)-s[n]\bar{h}_{n-1}(qx,s,q). \tag{2.16}$$

### 2.2. Generating functions

In the following we need different $q-$ analogues of the exponential series. These are well known, but for the convenience of the reader I state them explicitly.

The power series

$$e_q(z)=\sum_{n\geq 0}\frac{z^n}{[n]!} \tag{2.17}$$

and

$$E_q(z)=\sum_{n\geq 0}q^{\binom{n}{2}}\frac{z^n}{[n]!} \tag{2.18}$$

are related by

$$e_q(z)E_q(-z)=1 \tag{2.19}$$

and

$$e_{q^{-1}}(z)=E_q(z). \tag{2.20}$$

They can be regarded as formal power series or as convergent power series. In the second case (2.17) converges for $|z|<1$ whereas (2.18) is an entire function.

Let us note that

$$\begin{aligned}
&D_q e_q(ax)=ae_q(ax),\\
&D_q E_q(ax)=aE_q(qax),\\
&D_q e_{q^2}\left(\frac{ax^2}{1+q}\right)=axe_{q^2}\left(\frac{ax^2}{1+q}\right),\\
&D_q E_{q^2}\left(\frac{ax^2}{1+q}\right)=axE_{q^2}\left(\frac{q^2ax^2}{1+q}\right).
\end{aligned} \tag{2.21}$$

From the definition of $D_q$ these identities are equivalent with

$$\begin{aligned}
&e_q(x)\big(1-(1-q)x\big)=e_q(qx),\\
&E_q(x)=\big(1+(1-q)x\big)E_q(qx).
\end{aligned} \tag{2.22}$$

This implies the series expansion

$$e(z,q)=\sum_n\frac{z^n}{(q;q)_n}=\frac{1}{(z;q)_\infty}, \tag{2.23}$$

for $|z|<1$ and the expansion

$$E(z,q)=\sum_n q^{\binom{n}{2}}\frac{z^n}{(q;q)_n}=(-z;q)_\infty \tag{2.24}$$

which converges for all $z\in\mathbb{C}.$

In their domains of convergence we have $e_q(z)=e((1-q)z,q)$ and $E_q(z)=E((1-q)z,q).$

Note that

$$e(z,q^{-1})=\left(qz;q\right)_\infty, \tag{2.25}$$

because

$$e(z,q^{-1})=\sum_n\frac{z^n}{(q^{-1};q^{-1})_n}=\sum_n(-1)^n q^{\binom{n+1}{2}}\frac{z^n}{(q;q)_n}=E\left(-qz,q\right)=\left(qz;q\right)_\infty.$$

Let us also note that (cf. e.g. [4])

$$\frac{e_q(xz)}{e_q(az)}=\sum_{n\geq 0}\left(x-a\right)\left(x-qa\right)\cdots\left(x-q^{n-1}a\right)\frac{z^n}{[n]!}. \tag{2.26}$$

As a $q-$ analogue of (1.5) we get the generating function of the $q-$ Hermite polynomials

$$\sum_{n\geq 0}H_n(x,s,q)\frac{z^n}{[n]!}=\frac{e_q(xz)}{e_{q^2}\left(\frac{qsz^2}{1+q}\right)}=e_q(xz)E_{q^2}\left(\frac{-qsz^2}{1+q}\right). \tag{2.27}$$

For (2.4) is equivalent with

$$\sum_{n\geq 0}H_n(x,s,q)\frac{z^n}{[n]!}=\sum_j\frac{x^j z^j}{[j]!}\sum_k(-s)^k q^{k^2}[2k-1]!!\frac{z^{2k}}{[2k]!}$$

and

$$\sum_k(-s)^k q^{k^2}\frac{z^{2k}}{[2k]_q!}[2k-1]_q!!=\sum_k(-s)^k q^{k^2}\frac{z^{2k}}{[2]_q[4]_q\cdots[2k]_q}$$

$$=\sum_k(-1)^k q^{k(k-1)}\frac{(qs)^k}{(1+q)^k}\frac{z^{2k}}{[k]_{q^2}!}=\frac{1}{e_{q^2}\left(\frac{qsz^2}{1+q}\right)}.$$

An equivalent version of the generating function for $|z|<1$ is

$$\sum_{n\geq 0} H_n(x,(1-q)s,q)\frac{z^n}{(q;q)_n} = e(xz,q)E\left(-qsz^2,q^2\right) = \frac{\left(qsz^2;q^2\right)_\infty}{\left(xz;q\right)_\infty}. \tag{2.28}$$

Note that $\lim_{q\to 1} H_n(x,(1-q)s,q) = x^n$. Whereas this is no direct $q-$ analogue of the Hermite polynomials it is very useful since the infinite products on the right-hand side are sometimes easier to handle.

Since $D_q e_q(xz) = ze_q(xz)$ we see that (2.27) is equivalent with

$$H_n(x,s,q) = E_{q^2}\left(\frac{-qsD_q^2}{1+q}\right)x^n, \tag{2.29}$$

which generalizes (1.6).

From (2.27) we see that $e_q(xz) = e_{q^2}\left(\frac{qsz^2}{1+q}\right)\sum_{j\geq 0} H_j(x,s,q)\frac{z^j}{[j]_q!}$ and therefore

$$x^n = \sum_{k=0}^{\left\lfloor \frac{n}{2}\right\rfloor} (qs)^k \begin{bmatrix} n \\ 2k \end{bmatrix} [2k-1]!! H_{n-2k}(x,s,q). \tag{2.30}$$

Since $\frac{d}{dz}e^{zx-s\frac{z^2}{2}} = \left(x-sz\right)e^{zx-s\frac{z^2}{2}}$ we have

$$\sum_{n\geq 0} H_{n+1}(x,s)\frac{z^n}{n!} = (x-sz)\sum_{n\geq 0} H_n(x,s)\frac{z^n}{n!}. \tag{2.31}$$

For $H_n(x,s,q)$ we deduce from (2.27) the $q-$ analogue

$$\sum_{n\geq 0} H_{n+1}(x,s,q)\frac{z^n}{[n]!} = \frac{x-qsz}{1-(1-q)xz}\sum_{n\geq 0} H_n(x,s,q)\frac{z^n}{[n]!}. \tag{2.32}$$

To prove this let $\frac{\partial}{\partial z} f(z) = \frac{f(z)-f(qz)}{z}$. Applying $\frac{\partial}{\partial z}$ to (2.28) we get

$$\sum_{n\ge 0} H_{n+1}(x,(1-q)s,q)\frac{z^n}{(q;q)_n}=\frac{\partial}{\partial z}\sum_{n\ge 0} H_n(x,(1-q)s,q)\frac{z^n}{(q;q)_n}=\frac{\partial}{\partial z}\frac{\left(qsz^2;q^2\right)_\infty}{(xz;q)_\infty}$$

$$=\frac{1}{z}\left(\frac{\left(qsz^2;q^2\right)_\infty}{(xz;q)_\infty}-\frac{\left(q^3sz^2;q^2\right)_\infty}{(qxz;q)_\infty}\right)=\frac{1}{z}\frac{\left(q^3sz^2;q^2\right)_\infty}{(qxz;q)_\infty}\left(\frac{1-qsz^2}{1-xz}-1\right)$$

$$=\frac{x-qsz}{1-xz}\frac{\left(q^3sz^2;q^2\right)_\infty}{(qxz;q)_\infty}=\frac{x-qsz}{1-xz}\sum_{n\ge 0} H_n(x,(1-q)s,q)\frac{q^n z^n}{(q;q)_n}.$$

Replacing $z\to(1-q)z$ we get (2.32).

In the same way as above we get

$$\sum_{n\ge 0} K_n(x,s,q)\frac{z^n}{[n]!}=e_{q^2}\left(\frac{-sz^2}{1+q}\right)E_q(xz) \tag{2.33}$$

and

$$\sum_{n\ge 0} K_n(x,(1-q)s,q)\frac{z^n}{(q;q)_n}=e\left(-sz^2,q^2\right)E(xz,q)=\frac{(-xz;q)_\infty}{\left(-sz^2;q^2\right)_\infty}. \tag{2.34}$$

In this case

$$\sum_{n\ge 0} K_{n+1}(x,(1-q)s,q)\frac{z^n}{(q;q)_n}=\frac{\partial}{\partial z}\sum_{n\ge 0} K_n(x,(1-q)s,q)\frac{z^n}{(q;q)_n}=\frac{\partial}{\partial z}\frac{(-xz;q)_\infty}{\left(-sz^2;q^2\right)_\infty}$$

$$=\frac{1}{z}\left(\frac{(-xz;q)_\infty}{\left(-sz^2;q^2\right)_\infty}-\frac{(-qxz;q)_\infty}{\left(-sq^2z^2;q^2\right)_\infty}\right)=\frac{1}{z}\frac{(-qxz;q)_\infty}{\left(-sq^2z^2;q^2\right)_\infty}\left(\frac{1+xz}{1+sz^2}-1\right)$$

$$=\frac{x-sz}{1+sz^2}\frac{(-qxz;q)_\infty}{\left(-sq^2z^2;q^2\right)_\infty}=\frac{x-sz}{1+sz^2}\sum_{n\ge 0} K_n(x,(1-q)s,q)\frac{(qz)^n}{(q;q)_n}.$$

**Remark**

For $\overline{h}_n(x,s,q)$ we get

$$\sum_{n\ge 0} \overline{h}_n(x,s,q)\frac{z^n}{[n]!}=e_q(xz)e_{q^2}\left(\frac{-sz^2}{1+q}\right) \tag{2.35}$$

and thus

$$\sum_{n\ge 0} \overline{h}_n(x,(1-q)s,q)\frac{z^n}{(q;q)_n}=e(xz,q)e\left(-sz^2,q^2\right)=\frac{1}{(xz;q)_\infty\left(-sz^2;q^2\right)_\infty}. \tag{2.36}$$

From this we see that for the linear functional $\lambda$ defined by $\lambda\left(\overline{h}_n(x,s,q)\right)=[n=0]$ we get the moments

$\lambda\left(x^{2n+1}\right)=0$ and $\lambda\left(x^{2n}\right)=q^{n^2-n}s^n[2n-1]!!.$

It is easily seen that $\lambda\left(\overline{h}_1\overline{h}_3\right)=\left(q^2-1\right)[3]s^2\neq 0,$ which implies that the sequence $\left(\overline{h}_n\right)$ cannot be orthogonal.

### 2.3. The simplest special cases

**2.3.1.** The polynomials with the simplest right-hand side of (2.28) occur for $qs=1$.

Here we get the polynomials

$$h_n(x;q)=H_n\left(x,\frac{1-q}{q},q\right) \tag{2.37}$$

which have been called in [13] **discrete $q$ – Hermite polynomials I**.

They satisfy

$$h_{n+1}(x;q)=xh_n(x;q)-q^{n-1}\left(1-q^n\right)h_{n-1}(x;q) \tag{2.38}$$

with generating function for $|xz|<1$

$$\sum_{n\geq 0}h_n(x;q)\frac{z^n}{(q;q)_n}=\frac{e(xz,q)}{e\left(z^2,q^2\right)}=\frac{\left(z^2;q^2\right)_\infty}{(xz;q)_\infty}. \tag{2.39}$$

The polynomials $h_n(x;q)$ are given by

$$h_n(x;q)=\sum_{2k\leq n}(-1)^k q^{k^2-k}\begin{bmatrix}n\\2k\end{bmatrix}\left(q;q^2\right)_k x^{n-2k}. \tag{2.40}$$

The first terms are

$$\left(h_n(x;q)\right)_{n\geq 0}=\left(1,\ x,\ x^2-(1-q),\ x^3-\left(1-q^3\right)x,\ x^4-(1-q)\begin{bmatrix}4\\2\end{bmatrix}x^2+q^2(1-q)\left(1-q^3\right),\cdots\right)$$

These polynomials have first been considered by Al-Salam and Carlitz in [1]. They studied more generally polynomials $U_n^{(a)}(x)$ with generating function

$$\sum_{n\geq 0}U_n^{(a)}(x)\frac{z^n}{(q;q)_n}=\frac{e(xz,q)}{e(z,q)e(az,q)},$$

which for $a=-1$ reduce to $h_n(x;q)$, because $\frac{1}{e(z,q)e(-z,q)}=(z;q)_\infty(-z;q)_\infty=(z^2;q^2)_\infty$.

From

$$\sum_{n\geq 0} h_n(x;q)\frac{z^n}{(q;q)_n}=\frac{(z;q)_\infty(-z;q)_\infty}{(xz;q)_\infty}=(-z;q)_\infty\frac{(z;q)_\infty}{(xz;q)_\infty}$$
$$=\sum_j q^{\binom{j}{2}}\frac{z^j}{(q;q)_j}\sum_k (x-1)(x-q)\cdots(x-q^{k-1})\frac{z^k}{(q;q)_k}$$

we deduce by (2.26) that

$$h_n(x;q)=\sum_{j=0}^{n} q^{\binom{n-j}{2}}\begin{bmatrix} n \\ j \end{bmatrix}(x-1)(x-q)\cdots(x-q^{j-1}). \tag{2.41}$$

The identity (2.29) becomes

$$\frac{1}{e(\delta_q^2,q^2)}x^n=h_n(x,q) \tag{2.42}$$

where $\delta_q=(1-q)D_q$ with $\delta_q x^n=(1-q^n)x^{n-1}$.

By (2.30) we have

$$x^n=\sum_{k=0}^{\left\lfloor \frac{n}{2}\right\rfloor}\begin{bmatrix} n \\ 2k \end{bmatrix}(q;q^2)_k\, h_{n-2k}(x;q). \tag{2.43}$$

For the linear functional $\Lambda$ defined by $\Lambda(h_n(x;q))=[n=0]$ we get the moments

$$\begin{aligned}\Lambda(x^{2n})&=(q;q^2)_n\\ \Lambda(x^{2n+1})&=0.\end{aligned} \tag{2.44}$$

**2.3.2.** Instead of the polynomials $h_n(x;q)$ for $q>1$ it is convenient to consider the polynomials

$$\widehat{h}_n(x;q)=i^{-n}h_n\left(ix;\frac{1}{q}\right)=i^{-n}H_n\left(ix,q-1,\frac{1}{q}\right)=H_n\left(x,1-q,\frac{1}{q}\right)=q^{-\binom{n}{2}}K_n(x,1-q,q) \tag{2.45}$$

for $0<q<1$, the so called **discrete $q-$ Hermite polynomials II**.

From (2.11) we see that the polynomials $\widehat{h}_n(x;q)$ satisfy

$$\widehat{h}_{n+1}(x;q)=x\widehat{h}_n(x;q)-q^{-2n+1}(1-q^n)\widehat{h}_{n-1}(x;q). \tag{2.46}$$

By (2.34) we have

$$\sum_{n\ge 0}\widehat{h}_n(x;q)q^{\binom{n}{2}}\frac{z^n}{(q;q)_n}=\frac{(-xz;q)_\infty}{\left(-z^2;q^2\right)_\infty}. \tag{2.47}$$

These polynomials have the explicit expression

$$q^{\binom{n}{2}}\widehat{h}_n(x;q)=\sum_{2k\le n}(-1)^k q^{\binom{n-2k}{k}}\begin{bmatrix} n \\ 2k \end{bmatrix}\left(q;q^2\right)_k x^{n-2k}. \tag{2.48}$$

It is easily verified that

$$\sum_{2k\le n}q^{3k^2-2kn}\begin{bmatrix} n \\ 2k \end{bmatrix}\left(q;q^2\right)_k \widehat{h}_{n-2k}(x;q)=x^n. \tag{2.49}$$

Therefore the moments with respect to the linear functional $L$ defined by

$$L\left(\widehat{h}_n(x;q)\right)=[n=0] \tag{2.50}$$

are

$$L\left(x^{2n}\right)=\frac{\left(q;q^2\right)_n}{q^{n^2}} \qquad L\left(x^{2n+1}\right)=0. \tag{2.51}$$

**2.3.3.** An analog of (2.41) for $\overline{h}_n(x,q-1,q)$ is

$$\overline{h}_n(x,q-1,q)=\sum_{k=0}^{n}\begin{bmatrix} n \\ k \end{bmatrix}\sum_{j=0}^{k}\begin{bmatrix} k \\ j \end{bmatrix}(-1)^{k-j}x^j. \tag{3.1}$$

This immediately follows from (2.36).

## 3. Rodrigues-type formulae

### 3.1. Bivariate and discrete q- Hermite polynomials I

**3.1.1.** For the bivariate $q-$ Hermite polynomials we get two different $q-$ analogues of (1.10). The first one is

$$H_n(x,s,q)=\left(\mathbf{x}-q^{n-1}sD_q\right)\left(\mathbf{x}-q^{n-2}sD_q\right)\cdots\left(\mathbf{x}-sD_q\right)1. \tag{3.2}$$

This follows from (2.1) and (2.3) which imply

$$H_n(x,s,q)=xH_{n-1}(x,s,q)-q^{n-1}sD_qH_{n-1}(x,s,q)=\left(\mathbf{x}-q^{n-1}sD_q\right)H_{n-1}(x,s,q).$$

The second one is

$$H_n(x,s,q)=\left(\mathbf{x}-qsD_q\right)\left(\mathbf{x}-q^3sD_q\right)\cdots\left(\mathbf{x}-q^{2n-1}sD_q\right)1. \qquad (3.3)$$

To prove this we observe that

$H_{n+1}(x,s,q)=xH_n(x,s,q)-q^ns[n]H_{n-1}(x,s,q)$ and $\left(D_qH_{n+1}\right)(x,s,q)=[n+1]H_n(x,s,q)$

imply $H_{n+1}(qx,s,q)-H_{n+1}(x,s,q)=(q-1)x[n+1]H_n(x,s,q)=\left(q^{n+1}-1\right)xH_n(x,s,q)$

and therefore

$$H_{n+1}(qx,s,q)=xH_n(x,s,q)-q^n[n]sH_{n-1}(x,s,q)+\left(q^{n+1}-1\right)xH_n(x,s,q)$$
$$=q^{n+1}xH_n(x,s,q)-q^n[n]sH_{n-1}(x,s,q)$$

Changing $x\to\frac{x}{q}$ we get another recurrence relation for the $q-$ Hermite polynomials

$$H_{n+1}(x,s,q)=q^nxH_n\left(\frac{x}{q},s,q\right)-q^n[n]sH_{n-1}\left(\frac{x}{q},s,q\right). \qquad (3.4)$$

Since $H_n(x,q^2s,q)=q^nH_n\left(\frac{x}{q},s,q\right)$ this is equivalent with

$$H_{n+1}(x,s,q)=xH_n\left(x,q^2s,q\right)-qs[n]H_{n-1}\left(x,q^2s,q\right). \qquad (3.5)$$

This implies

$$H_n(x,s,q)=xH_{n-1}(x,q^2s,q)-qsD_qH_{n-1}(x,q^2s,q)=\left(\mathbf{x}-qsD_q\right)H_{n-1}(x,q^2s,q)$$
$$=\left(\mathbf{x}-qsD_q\right)\left(\mathbf{x}-q^3sD_q\right)\cdots\left(\mathbf{x}-q^{2n-1}sD_q\right)1.$$

**(3.4)** also implies that

$$q^{-\binom{n}{2}}H_n(x,s,q)=\left(\left(x-qsD_q\right)\varepsilon^{-1}\right)^n1. \qquad (3.6)$$

Let us give another formula for the right-hand side of this identity.

The $q-$ differentiation operator $D_q$ satisfies

$$D_q\left(f(x)g(x)\right)=f(qx)D_qg(x)+\left(D_qf\right)(x)g(x) \qquad (3.7)$$

and thus

$$D_q\mathbf{f(x)}=\mathbf{f(qx)}D_q+\left(\mathbf{D_qf}\right)\mathbf{(x)} \qquad (3.8)$$

and

$$D_q \mathbf{f(x)} = \mathbf{f(x)} D_q + \left(\mathbf{D_q f}\right)\mathbf{(x)}\varepsilon. \tag{3.9}$$

Since $D_q E_{q^2}\left(\frac{ax^2}{1+q}\right) = axE_{q^2}\left(\frac{q^2ax^2}{1+q}\right)$ we get from (3.8)

$$\frac{1}{\mathbf{E_{q^2}}\left(\frac{\mathbf{-aqx^2}}{\mathbf{1+q}}\right)}\varepsilon^{-1}(-sD_q)\mathbf{E_{q^2}}\left(\frac{\mathbf{-aqx^2}}{\mathbf{1+q}}\right) = -qsD_q\varepsilon^{-1} + as\mathbf{x}\varepsilon^{-1} \tag{3.10}$$

and therefore for $a = \frac{1}{s}$

$$q^{-\binom{n}{2}} H_n(x,s,q) = \frac{1}{E_{q^2}\left(\frac{-qx^2}{(1+q)s}\right)}(-qsD_q\varepsilon^{-1})^n E_{q^2}\left(\frac{-qx^2}{(1+q)s}\right). \tag{3.11}$$

Letting in (3.11) $s \to (1-q)s$ and $a \to \frac{1}{(1-q)s}$ we get

$$\frac{1}{E_{q^2}\left(\frac{-qx^2}{\left(1-q^2\right)s}\right)}\left(-(1-q)qsD_q\varepsilon^{-1}\right)^n E_{q^2}\left(\frac{-qx^2}{\left(1-q^2\right)s}\right) = q^{-\binom{n}{2}} H_n(x,(1-q)s,q)$$

or equivalently

$$(q-1)^n s^n \frac{1}{\left(\frac{qx^2}{s};q^2\right)_\infty}\left(\varepsilon^{-1}D_q\right)^n\left(\frac{qx^2}{s};q^2\right)_\infty = q^{-\binom{n}{2}} H_n(x,(1-q)s,q). \tag{3.12}$$

**3.1.2.**

By (3.12) we have

$$h_n(x;q) = q^{\binom{n}{2}-n}(q-1)^n \frac{1}{\left(q^2x^2;q^2\right)_\infty}\left(\varepsilon^{-1}D_q\right)^n\left(q^2x^2;q^2\right)_\infty. \tag{3.13}$$

Let us give another proof.

Since $(1-q)D_q\left(q^2x^2;q^2\right)_\infty = -q^2x\left(q^4x^2;q^2\right)_\infty$

we get by (3.7)

$$(1-q)\varepsilon^{-1}D_q\left(\mathbf{q^2x^2;q^2}\right)_\infty=-qx\left(\mathbf{q}^2\mathbf{x^2;q^2}\right)_\infty\varepsilon^{-1}+\left(\mathbf{q}^2\mathbf{x^2;q^2}\right)_\infty(1-q)qD_q\varepsilon^{-1}\ .$$

or

$$\frac{1}{\left(\mathbf{q}^2\mathbf{x^2;q^2}\right)_\infty}\frac{q-1}{q}\varepsilon^{-1}D_q\left(\mathbf{q}^2\mathbf{x^2;q^2}\right)_\infty=\left(\mathbf{x}\varepsilon^{-1}-(1-q)D_q\varepsilon^{-1}\right). \tag{3.14}$$

This implies (3.13)

Observing that $D_{q^{-1}}=\varepsilon^{-1}D_q$ this formula implies [13], (3.28.9)

$$D_{q^{-1}}\left(q^2x^2;q^2\right)_\infty h_n(x;q)=-\frac{q^{-n+1}}{1-q}\left(q^2x^2;q^2\right)_\infty h_{n+1}(x;q). \tag{3.15}$$

### 3.2. The discrete q-Hermite polynomials II

By (3.6) for $q\to\frac{1}{q}$ we deduce the following $q-$ analogue of (1.10):

$$K_n(x,s,q)=\left(\mathbf{x}\varepsilon-sD_q\right)^n 1. \tag{3.16}$$

For $f(x)=e_{q^2}\left(\frac{ax^2}{[2]_q}\right)$ we have $D_q\mathbf{f(x)}=\mathbf{f(x)}D_q+a\mathbf{x}\varepsilon$ or

$$\mathbf{e}_{\mathbf{q}^2}\left(\frac{\mathbf{ax}^2}{[\mathbf{2}]_\mathbf{q}}\right)^{-1}D_q\mathbf{e}_{\mathbf{q}^2}\left(\frac{\mathbf{ax}^2}{[\mathbf{2}]_\mathbf{q}}\right)=D_q+a\mathbf{x}\varepsilon. \tag{3.17}$$

Comparing (3.17) with (3.16) we get the Rodrigues-type formula

$$K_n(x,s,q)=\left(\mathbf{x}\varepsilon-sD_q\right)^n 1=(-s)^n\mathbf{e}_{\mathbf{q}^2}\left(\frac{\mathbf{-x}^2}{[\mathbf{2}]_\mathbf{q}\mathbf{s}}\right)^{-1}D_q^n\mathbf{e}_{\mathbf{q}^2}\left(\frac{\mathbf{-x}^2}{[\mathbf{2}]_\mathbf{q}\mathbf{s}}\right). \tag{3.18}$$

Since $-sD_q\dfrac{1}{\left(-\frac{1-q}{s}x^2;q^2\right)_\infty}=-x\dfrac{1}{\left(-\frac{1-q}{s}x^2;q^2\right)_\infty}$

we get in the same way

$$\left(-\frac{\mathbf{1-q}}{\mathbf{s}}\mathbf{x}^2;\mathbf{q}^2\right)_\infty (-sD_q)\frac{\mathbf{1}}{\left(-\frac{\mathbf{1-q}}{\mathbf{s}}\mathbf{x}^2;\mathbf{q}^2\right)_\infty} = -sD_q + \mathbf{x}\varepsilon$$

and therefore the equivalent formula

$$K_n(x,s,q) = \left(\mathbf{x}\varepsilon - sD_q\right)^n 1 = (-s)^n \left((q-1)s^{-1}x^2;q^2\right)_\infty D_q^n \frac{1}{\left((q-1)s^{-1}x^2;q^2\right)_\infty}. \tag{3.19}$$

For $s = 1-q$ this reduces to the Rodrigues type formula [13], (3.29.10)

$$\frac{1}{\left(-x^2;q^2\right)_\infty}\widehat{h}_n(x;q) = (q-1)^n q^{-\binom{n}{2}} D_q^n \frac{1}{\left(-x^2;q^2\right)_\infty}. \tag{3.20}$$

**Remark**

For the polynomials $\overline{h}_n(x,s,q)$ we get (cf. [12]) by (2.14) and (2.16).

$$\overline{h}_n(x,s,q) = \left(x - s\varepsilon D_q\right)^n 1. \tag{3.21}$$

The sequence of polynomials $\left(\overline{h}_n(x,-qs,q)\right)$ is the umbral-inverse sequence to the sequence $\left(H_n(x,s,q)\right)$. This means that the linear map $U$ defined by $U\left(x^n\right) = \overline{h}_n(x,s,q) = \left(x - s\varepsilon D_q\right)^n 1$ satifies $U\left(H_n(x,s,q)\right) = x^n$ or with other words that

$$H_n\left(\mathbf{x} + qs\varepsilon D_q\right)1 = x^n. \tag{3.22}$$

The proof is straightforward by induction since (3.22) hold for $n = 0$ and $n = 1$:

$H_{n+1}\left(\mathbf{x} + qs\varepsilon D_q, s, q\right) = \left(\mathbf{x} + qs\varepsilon D_q\right)x^n - q^n[n]sx^{n-1} = x^{n+1}.$

In another form this fact has already been shown in (2.30).

The Hermite polynomials satisfy (cf. e.g. [15])

$$\sum_{k=0}^{n}\binom{n}{k}H_{n-k}(x,s)H_k(y,-s) = (x+y)^n. \tag{3.23}$$

A $q-$ analogue is

$$\sum_{k=0}^{n}\begin{bmatrix} n \\ k \end{bmatrix} x^k y^{n-k} = \sum_{k=0}^{n}\begin{bmatrix} n \\ k \end{bmatrix} H_{n-k}(x,s,q)\overline{h}_k(y,-qs,q). \tag{3.24}$$

This follows from

$$e_q(xz)e_q(yz)=\frac{e_q(xz)}{e_{q^2}\left(\frac{qsz^2}{1+q}\right)}e_q(yz)e_{q^2}\left(\frac{qsz^2}{1+q}\right)=\sum_{n\geq 0}H_n(x,s,q)\frac{z^n}{[n]!}\sum_{n\geq 0}\overline{h}_n(y,-qs,q)\frac{z^n}{[n]!}$$

by comparing coefficients.

### 3.3. A q-Burchnall formula

In [6] I have given some $q-$ analogues of Burchnall's formula. The simplest one is the following $q-$ Burchnall formula

$$\left(\mathbf{x}-qsD_q\right)\left(\mathbf{x}-q^3sD_q\right)\cdots\left(\mathbf{x}-q^{2n-1}sD_q\right)=\sum_{k=0}^{n}\begin{bmatrix}n\\k\end{bmatrix}q^{kn}H_{n-k}(\mathbf{x},s,q)\left(-sD_q\right)^k. \qquad (3.25)$$

The proof is by induction.

$$\left(\mathbf{x}-qsD_q\right)\sum_{k=0}^{n}\begin{bmatrix}n\\k\end{bmatrix}q^{kn}H_{n-k}(\mathbf{x},q^2s,q)\left(-q^2sD_q\right)^k=\sum_{k=0}^{n}\begin{bmatrix}n\\k\end{bmatrix}q^{kn}\mathbf{x}H_{n-k}(\mathbf{x},q^2s,q)\left(-q^2sD_q\right)^k$$
$$-qs\sum_{k=0}^{n}\begin{bmatrix}n\\k\end{bmatrix}q^{k(n+1)}H_{n-k}(q\mathbf{x},q^2s,q)\left(-qsD_q\right)^kD_q-qs\sum_{k=0}^{n}\begin{bmatrix}n\\k\end{bmatrix}q^{kn}\left(D_qH_{n-k}(\mathbf{x},q^2s,q)\right)\left(-q^2sD_q\right)^k$$
$$=\sum_{k=0}^{n}\begin{bmatrix}n\\k\end{bmatrix}q^{kn}\left(\left(\mathbf{x}-qsD_q\right)H_{n-k}(\mathbf{x},q^2s,q)\right)\left(-q^2sD_q\right)^k+\sum_{k}\begin{bmatrix}n\\k\end{bmatrix}q^{k(n+1)}H_{n-k}(q\mathbf{x},q^2s,q)\left(-qsD_q\right)^{k+1}$$
$$=\sum_{k}q^k\begin{bmatrix}n\\k\end{bmatrix}q^{k(n+1)}H_{n+1-k}(\mathbf{x},s,q)\left(-sD_q\right)^k+\sum_{k}\begin{bmatrix}n\\k\end{bmatrix}q^{(k+1)(n+1)}H_{n-k}(\mathbf{x},s,q)\left(-sD_q\right)^{k+1}$$
$$=\sum_{k}\left(q^k\begin{bmatrix}n\\k\end{bmatrix}+\begin{bmatrix}n\\k-1\end{bmatrix}\right)q^{k(n+1)}H_{n+1-k}(\mathbf{x},s,q)\left(-sD_q\right)^k=\sum_{k}\begin{bmatrix}n+1\\k\end{bmatrix}q^{k(n+1)}H_{n+1-k}(\mathbf{x},s,q)\left(-sD_q\right)^k.$$

If we apply (3.25) to $H_m\left(x,q^{2n}s,q\right)$ we get a $q-$ analogue of Nielsen's identity

$$H_{n+m}(x,s,q)=\sum_{k}\begin{bmatrix}n\\k\end{bmatrix}\begin{bmatrix}m\\k\end{bmatrix}[k]!q^{kn}(-s)^kH_{n-k}(x,s,q)H_{m-k}\left(x,q^{2n}s,q\right). \qquad (3.26)$$

Since $H_{m-k}\left(x,q^{2n}s,q\right)=q^{(m-k)n}H_{m-k}\left(\frac{x}{q^n},s,q\right)$ this coincides with

$$H_{n+m}(x,s,q)=q^{mn}\sum_{k}\begin{bmatrix}n\\k\end{bmatrix}\begin{bmatrix}m\\k\end{bmatrix}[k]!(-s)^kH_{n-k}(x,s,q)H_{m-k}\left(\frac{x}{q^n},s,q\right), \qquad (3.27)$$

which has been proved in [4] with another method.

## 4. Associated probability measures

### 4.1. The q-integral

In the classical case the linear functional $\Lambda$ which satisfies $\Lambda\left(H_n(x,s)\right)=[n=0]$ is given by

$$\Lambda\left(f\right)=\frac{1}{\sqrt{2\pi s}}\int_{-\infty}^{\infty} f(x)e^{-\frac{x^2}{2s}}dx.$$

In order to find a $q-$ analogue we need the Jackson $q-$ integral (cf. [10], [11], [14]). We assume that $0<q<1.$

We call $F(x)$ a $q-$ antiderivative of $f(x)$ if $D_qF(x)=f(x).$ This means that
$F(x)-F(qx)=(1-q)xf(x)$ or $(1-\varepsilon_q)F(x)=(1-q)xf(x)$ or
$F(x)=\frac{1}{1-\varepsilon_q}(1-q)xf(x)=\sum_{n\geq 0}\varepsilon_q^n(1-q)xf(x)=\sum_{n\geq 0}(1-q)q^nxf(q^nx).$

If this sum converges absolutely it is clear that $F(x)$ is a $q-$ antiderivative of $f(x)$ because $F(x)-F(qx)=(1-q)xf(x).$ In the classical case all antiderivatives of $0$ are constants. In the $q-$ case also each function $\varphi(x)$ with $\varphi(qx)=\varphi(x)$ is a $q-$ antiderivative of $0.$ But it can be shown that up to a constant any function has at most one $q-$ antiderivative that is continuous at $x=0.$ (Cf. [11] for details).

This leads to the following definition. We always assume that $0<q<1.$

Let $0<a<b.$ The definite $q-$integral is defined as

$$\int_0^b f(x)d_qx=\sum_{j=0}^{\infty} f\left(q^jb\right)\left(q^j-q^{j+1}\right)b \tag{4.1}$$

and

$$\int_a^b f(x)d_qx=\int_0^b f(x)d_qx-\int_0^a f(x)d_qx \tag{4.2}$$

provided that the sums converge absolutely, for example if $\left|x^\alpha f(x)\right|\leq M$ in a neighbourhood of $0$ for some $0\leq\alpha<1.$ For then we have $\left|q^j f(q^jx)\right|<Mq^j\left(q^jx\right)^{-\alpha}=Mx^{-\alpha}\left(q^{1-\alpha}\right)^j.$

Note that $\int_a^b f(x)d_qx$ depends on the values of $f(x)$ in the whole interval $(0,b].$

If $f(x)$ is continuous at $x=0$ and the $q-$integral converges absolutely, then

$$\int_a^b D_qf(x)d_qx=f(b)-f(a). \tag{4.3}$$

For all $q-$antiderivatives of $D_qf(x)$ which are continuous at $0$ are of the form $f(x)+C.$

From the product rule

$$D_q\left(f(x)g(x)\right)=f(x)\left(D_q g(x)\right)+g(qx)\left(D_q f(x)\right)$$

we obtain by (4.3) the formula for integration by parts

$$\int_a^b f(x)\left(D_q g(x)\right)d_q x=f(b)g(b)-f(a)g(a)-\int_a^b g(qx)\left(D_q f(x)\right)d_q x. \qquad (4.4)$$

## 4.2. The probability measure for the discrete q-Hermite polynomials I

We start from the Rodrigues-type formula

$$h_n(x;q)=q^{\binom{n}{2}-n}(q-1)^n\frac{1}{\left(q^2x^2;q^2\right)_\infty}\left(\varepsilon^{-1}D_q\right)^n\left(q^2x^2;q^2\right)_\infty.$$

By multiplying both sides with $\left(q^2x^2;q^2\right)_\infty$ and applying $\varepsilon^{-1}$ we see that

$\left(\varepsilon^{-1}D_q\right)^n\left(x^2;q^2\right)_\infty$ is of the form $C_n h_n\left(\frac{x}{q};q\right)\left(x^2;q^2\right)_\infty$ and thus vanishes at $x=\pm 1$.

This implies that for $n>0$

$$\int_{-1}^1 h_n(x;q)\left(q^2x^2;q^2\right)_\infty d_q x=0. \qquad (4.5)$$

Since $\left(q^2x^2;q^2\right)_\infty$ is continuous and therefore bounded on $[-1,1]$ the $q-$integral (4.5) is given by an absolutely convergent series.

Thus the linear functional $\Lambda$ which gives $\Lambda\left(h_n(x;q)\right)=[n=0]$ is

$$\Lambda\left(f(x)\right)=\frac{\int_{-1}^1 f(x)\left(q^2x^2;q^2\right)_\infty d_q x}{\int_{-1}^1\left(q^2x^2;q^2\right)_\infty d_q x}. \qquad (4.6)$$

Let us now calculate $\int_{-1}^1\left(q^2x^2;q^2\right)_\infty d_q x=2\int_0^1\left(q^2x^2;q^2\right)_\infty d_q x.$

By definition of the $q-$ integral we have

$$\int_0^1\left(q^2x^2;q^2\right)_\infty d_q x=(1-q)\sum_n q^n\left(q^{2n};q^2\right)_\infty=(1-q)\left(q^2;q^2\right)_\infty\sum_n\frac{q^n}{\left(q^2;q^2\right)_n}$$

$$=(1-q)\left(q^2;q^2\right)_\infty e\left(q,q^2\right)=(1-q)\left(q;q\right)_\infty(-q;q)_\infty\frac{1}{\left(q;q^2\right)_\infty}.$$

As Euler has shown

$$\frac{1}{\left(q;q^2\right)_\infty}=\frac{1}{(1-q)\left(1-q^3\right)\left(1-q^5\right)\cdots}=\frac{(1-q^2)\left(1-q^4\right)\left(1-q^6\right)\cdots}{(1-q)\left(1-q^2\right)\left(1-q^3\right)\cdots}=\left(1+q\right)\left(1+q^2\right)\left(1+q^3\right)\cdots.$$

Therefore

$$\int_0^1\left(q^2x^2;q^2\right)_\infty d_qx=(1-q)\left(q;q\right)_\infty(-q;q)_\infty\left(-q;q\right)_\infty.$$

If we choose $x=-q$ in Jacobi's triple product identity

$$\sum_{k\in\mathbb{Z}}(-1)^k q^{\binom{k}{2}}x^k=(x;q)_\infty\left(\frac{q}{x};q\right)_\infty(q;q)_\infty$$

we see that

$$\int_{-1}^1\left(q^2x^2;q^2\right)_\infty d_qx=2(1-q)\left(q^2;q^2\right)_\infty\left(-q;q\right)_\infty=(1-q)\sum_{n=0}^\infty q^{\binom{n+1}{2}}. \tag{4.7}$$

Thus we get

$$\Lambda\left(f(x)\right)=\frac{\int_{-1}^1 f(x)\left(q^2x^2;q^2\right)_\infty d_qx}{\int_{-1}^1\left(q^2x^2;q^2\right)_\infty d_qx}=\frac{1}{\sum_{n=0}^\infty q^{\binom{n+1}{2}}}\sum_{j=0}^\infty\left(f(q^j)+f(-q^j)\right)q^j\left(q^{2+2j};q^2\right)_\infty. \tag{4.8}$$

Let us look what this gives for the moments $\Lambda\left(x^m\right)$. For odd $m$ it is clear that $\Lambda\left(x^m\right)=0$.

Observe that

$$\sum_{j=0}^\infty q^{(2m+1)j}\left(q^{2+2j};q^2\right)_\infty=\left(q^2;q^2\right)_\infty\sum_n\frac{q^{2(m+1)j}}{\left(q^2;q^2\right)_n}=\left(q^2;q^2\right)_\infty e\left(q^{2m+1};q^2\right)=\frac{\left(q^2;q^2\right)_\infty}{\left(q^{2m+1};q^2\right)_\infty}$$
$$=\frac{\left(q^2;q^2\right)_\infty}{\left(q;q^2\right)_\infty}\left(q;q^2\right)_m=\left(q^2;q^2\right)_\infty(-q;q)_\infty\left(q;q^2\right)_m=\left(q;q^2\right)_m\sum_{n=0}^\infty q^{\binom{n+1}{2}}.$$

This implies

$$\Lambda\left(x^{2m}\right)=\frac{1}{\sum_{n=0}^\infty q^{\binom{n+1}{2}}}\sum_{j=0}^\infty q^{(2m+1)j}\left(q^{2+2j};q^2\right)_\infty=\left(q;q^2\right)_m$$

which agrees with (2.44).

**Remark**

Formula (4.8) has already been found by Al-Salam and Carlitz [1].

They observed that the moments $\Lambda\left(x^m\right)$ satisfy

$\Lambda\left(x^{2m}\right)=\left(q;q^2\right)_m=(1-q^{2m-1})\Lambda\left(x^{2m-2}\right)$ and $\Lambda\left(x^{2m+1}\right)=0$

and therefore are characterized by $\Lambda\left(x^{m+1}\right)=\left(1-q^m\right)\Lambda\left(x^{m-1}\right)$ with initial values $\Lambda(1)=1$ and $\Lambda(x)=0$ and constructed an infinite sum with the same properties. This is simply the reverse of the computation above.

### 4.3. Probability measures for the discrete q-Hermite polynomials II

In this case we need the notion of an improper $q-$ integral.

We define for $c>0$

$$\int_0^{\infty.c} f(x)d_q(x)=\lim_{N\to\infty}\int_0^{q^{-N}c} f(x)d_q(x)=\lim_{N\to\infty}\sum_{j=0}^{\infty} f\left(q^{j-N}c\right)\left(q^j-q^{j+1}\right)q^{-N}c$$

$$=\lim_{N\to\infty}(1-q)\sum_{j=-N}^{\infty} f\left(q^jc\right)q^jc=(1-q)\sum_{j=-\infty}^{\infty} f\left(q^jc\right)q^jc$$

provided the last sum converges absolutely.

The formula

$$\int_a^b f(x)\left(D_qg(x)\right)d_qx=f(b)g(b)-f(a)g(a)-\int_a^b g(qx)\left(D_qf(x)\right)d_qx \tag{4.9}$$

also holds in the improper case by letting $b=q^{-N}c$ tending to infinity.

In a similar way we define

$$\int_{-\infty.c}^{\infty.c} f(x)d_q(x)=(1-q)\sum_{j=-\infty}^{\infty}\left(f\left(q^jc\right)+f\left(-q^jc\right)\right)q^jc. \tag{4.10}$$

By (3.20) we see that

$$\frac{1}{\left(-x^2;q^2\right)_\infty}\widehat{h}_n(x;q)=(q-1)^nq^{-\binom{n}{2}}D_q^n\frac{1}{\left(-x^2;q^2\right)_\infty}. \tag{4.11}$$

This implies that for $n>0$

$$\int_{-\infty.c}^{\infty.c}\frac{\widehat{h}_n(x;q)}{\left((-x^2;q^2\right)_\infty}d_qx=(q-1)^n q^{-\binom{n}{2}}\int_{-\infty.c}^{\infty.c}D_q^n\frac{1}{\left(-x^2;q^2\right)_\infty}=(q-1)^n D_q^{n-1}\frac{1}{\left(x^2;q^2\right)_\infty}\Bigg|_{-\infty.c}^{\infty.c}=0$$

if the integral exists and the limits are $0$.

For this it suffices to show that $\sum_{j=-\infty}^{\infty} q^j cf_m(q^j c)$ converges if $f_m(x)=\dfrac{x^m}{\left(-x^2;q^2\right)_\infty}$, $m\in\mathbb{N}$.

It is clear that $\sum_{j=0}^{\infty} q^j cf_m(q^j \mathrm{c})$ converges because $f_m(x)\le x^m$. Since

$\left(-x^2;q^2\right)_\infty \ge \left(q^{2m}x^2\right)^{m+1}=q^{2m^2+2m}x^{2m+2}$ we have

$$q^{-N}cf_m\left(q^{-N}c\right)\le q^{-N}c\frac{q^{-mN}c^m}{q^{2m^2+2m}q^{-2mN-2N}}=\frac{q^{(m+1)N}}{q^{2m^2+2m}}c^{m+1}$$

and therefore we see that also $\sum_{j=0}^{\infty} q^{-j} cf_m(q^{-j}\mathrm{c})$ converges.

Therefore in this case the linear functional $L$ has infinitely many representations as a $q-$ integral.

In order to compute

$$\int_{-\infty.c}^{\infty.c}\frac{1}{\left(-x^2;q^2\right)_\infty}d_qx=2(1-q)\sum_{j=-\infty}^{\infty}\frac{q^jc}{\left(-q^{2j}c^2;q^2\right)_\infty}=\frac{2(1-q)}{\left(-c^2;q^2\right)_\infty}\sum_{j=-\infty}^{\infty}q^jc\left(-c^2;q^2\right)_j$$

we need Ramanujan's summation formula (cf. [2], 10.5.3)

$$\sum_{n=-\infty}^{\infty}(a;q)_n x^n=\frac{(ax;q)_\infty\left(\frac{q}{ax};q\right)_\infty(q;q)_\infty}{(x;q)_\infty\left(\frac{q}{a};q\right)_\infty}.$$

This gives

$$\int_{-\infty.c}^{\infty.c}\frac{1}{\left((-x^2;q^2\right)_\infty}d_qx=\frac{2(1-q)}{\left(-c^2;q^2\right)_\infty}\sum_{j=-\infty}^{\infty}q^j\left(-c^2;q^2\right)_j=\frac{2(1-q)\left(-c^2q;q^2\right)_\infty\left(-\frac{q}{c^2};q^2\right)_\infty\left(q^2;q^2\right)_\infty}{\left(-c^2;q^2\right)_\infty\left(q;q^2\right)_\infty\left(-\frac{q^2}{c^2};q^2\right)_\infty}$$

Combining these results we see that the linear functional $L$ defined by $L\left(h_n(x;q)\right)=[n=0]$ can be represented by each of the improper $q-$ integrals (cf. [13])

$$L\left(f(x)\right)=\frac{\left(-c^2;q^2\right)_\infty\left(q;q^2\right)_\infty\left(-\frac{q^2}{c^2};q^2\right)_\infty}{2(1-q)\left(-c^2q;q^2\right)_\infty\left(-\frac{q}{c^2};q^2\right)_\infty\left(q^2;q^2\right)_\infty}\int_{-\infty.c}^{\infty.c}\frac{f(x)}{\left(-x^2;q^2\right)_\infty}d_qx. \tag{4.12}$$

As above we can compute

$$L\left(x^{2m}\right)=\frac{\left(q;q^2\right)_m}{q^{m^2}}, \tag{4.13}$$

which can also be obtained from (2.49).

## 5. Some curious identities

The tangent numbers $T_{2n+1}$ and Euler numbers $E_{2n}$ can be defined by

$$\frac{e^z-e^{-z}}{e^z+e^{-z}}=\sum_{n\geq 0}(-1)^n T_{2n+1}\frac{z^{2n+1}}{(2n+1)!}$$

and

$$\frac{2}{e^z+e^{-z}}=\sum_{n\geq 0}(-1)^n E_{2n}\frac{z^{2n}}{(2n)!}.$$

The first terms are

$\left(T_{2n+1}\right)_{n\geq 0}=\left(1,2,16,272,7936,\cdots\right)$ and $\left(E_{2n}\right)_{n\geq 0}=\left(1,1,5,61,1385,\cdots\right)$.

The polynomials $H_n(1,s)$ are polynomials in $s$ with degree $\left\lfloor \frac{n}{2}\right\rfloor$. Consider the linear functional $F$ on the polynomials in $s$ defined by

$$F\left(H_{2n}(1,s)\right)=[n=0]. \tag{5.1}$$

Then

$$F\left(H_{2n+1}(1,s)\right)=(-1)^n T_{2n+1} \tag{5.2}$$

and

$$F\left(s^n\right)=\frac{E_{2n}}{(2n-1)!!}. \tag{5.3}$$

To prove this let $H(z)=\sum_{n\geq 0}H_n(1,s)\frac{z^n}{n!}=e^{z-\frac{sz^2}{2}}$.

Then $H(z)=e^{2z}H(-z)$ and therefore

$$\sum_{n\geq 0}H_{2n+1}(1,s)\frac{z^{2n+1}}{(2n+1)!}=\frac{e^z-e^{-z}}{e^z+e^{-z}}\sum_{n\geq 0}H_{2n}(1,s)\frac{z^{2n}}{(2n)!}.$$

Applying $F$ we get

$$\sum_{n\geq 0}F\left(H_{2n+1}(1,s)\right)\frac{z^{2n+1}}{(2n+1)!}=\frac{e^z-e^{-z}}{e^z+e^{-z}}=\sum_{n\geq 0}(-1)^n T_{2n+1}\frac{z^{2n+1}}{(2n+1)!},$$

which gives (5.2).

To prove (5.3) observe that

$$F\left(e^{-\frac{sz^2}{2}}\right)=F\left(e^{-z}e^{z-\frac{sz^2}{2}}\right)=e^{-z}\sum_{n\geq 0}F\left(H_n(1,s)\right)\frac{z^n}{n!}=e^{-z}\left(1+\sum_{n\geq 0}F\left(H_{2n+1}(1,s)\right)\frac{z^{2n+1}}{(2n+1)!}\right)$$

$$=e^{-z}\left(1+\sum_{n\geq 0}(-1)^n T_{2n+1}\frac{z^{2n+1}}{(2n+1)!}\right)=e^{-z}\left(1+\frac{e^z-e^{-z}}{e^z+e^{-z}}\right)=\frac{2}{e^z+e^{-z}}=\sum_{n\geq 0}(-1)^n E_{2n}\frac{z^{2n}}{(2n)!}.$$

These arguments can immediately be transferred to the polynomials $H_n(1,s,q)$.

Define the $q-$ tangent and $q-$ Euler numbers by

$$\frac{e_q(z)-e_q(-z)}{e_q(z)+e_q(-z)}=\sum_{n\geq 0}(-1)^n T_{2n+1}(q)\frac{z^{2n+1}}{[2n+1]!}$$

and

$$\frac{2}{e_q(z)+e_q(-z)}=\sum_{n\geq 0}(-1)^n E_{2n}(q)\frac{z^{2n}}{[2n]!}.$$

Let now $\Phi$ be the linear functional defined by

$$\Phi\left(H_{2n}(1,s,q)\right)=[n=0]. \tag{5.4}$$

Then

$$\Phi\left(H_{2n+1}\right)=(-1)^n T_{2n+1}(q) \tag{5.5}$$

and

$$\Phi\left(s^n\right)=\frac{E_{2n}(q)}{q^{n^2}[2n-1]!!}. \tag{5.6}$$